\documentclass[12pt]{article}
\usepackage{amssymb}
\usepackage[small,nohug,heads=littlevee]{diagrams}
\usepackage{amsmath}
\diagramstyle[labelstyle=\scriptstyle]
\newcommand{\be}{\begin{equation}}
\newcommand{\ee}{\end{equation}}
\newcommand{\bea}{\begin{eqnarray}}
\newcommand{\eea}{\end{eqnarray}}
\newcommand{\ba}{\begin{array}}
\newcommand{\ea}{\end{array}}

\newcommand{\bc}{\begin{center}}
\newcommand{\ec}{\end{center}}
\newcommand{\ben}{\begin{enumerate}}
\newcommand{\een}{\end{enumerate}}
\newcommand{\bfi}{\begin{figure}}
\newcommand{\efi}{\end{figure}}

\newcommand{\bq}{\begin{quote}}
\newcommand{\eq}{\end{quote}}
\newcommand{\bqu}{\begin{quotation}}
\newcommand{\equ}{\end{quotation}}
\newenvironment{emphit}{\begin{itemize}}{\end{itemize}}
\newcommand{\bemp}{\begin{emphit}}
\newcommand{\eemp}{\end{emphit}}

\newcommand{\bt}{\begin{tabular}}
\newcommand{\et}{\end{tabular}}

\newtheorem{myth}{Theorem}[section]
\newtheorem{mylem}{Lemma}[section]
\newtheorem{mycor}{Corollary}[section]

\newtheorem{mydef}{Definition}[section]

\begin{document}
\date{}
\title{Holomorphy of Osborn loops
\footnote{2010 mathematics subject classification
primary 20N05; secondary 08A05.}\thanks{{\bf keywords: Osborn loops, holomorphy}}}
\author{A. O. Isere\\
Department of Mathematics,\\
Ambrose Alli University,
Ekpoma, Nigeria\\
abednis@yahoo.co.uk \and
J. O. Ad\'en\'iran \\
Department of Mathematics,\\
Federal University of Agriculture, \\
Abeokuta 110101, Nigeria.\\
ekenedilichineke@yahoo.com\\
adeniranoj@unaab.edu.ng \and
T. G. Jaiy\'e\d ol\'a\thanks{All correspondence to be addressed to this author.} \\
Department of Mathematics,\\
Obafemi Awolowo University,\\
Ile Ife 220005, Nigeria.\\
jaiyeolatemitope@yahoo.com\\tjayeola@oauife.edu.ng}\maketitle
\begin{abstract}
Let $(L,\cdot)$ be any loop and let $A(L)$ be a group of automorphisms of $(L,\cdot)$ such that $\alpha$ and $\phi$ are elements of $A(L)$. It is shown that, for all
$x,y,z\in L$, the $A(L)$-holomorph $(H,\circ)=H(L)$ of $(L,\cdot)$ is an Osborn loop if and only if
$x\alpha (yz\cdot x\phi^{-1})= x\alpha (yx^\lambda\cdot x) \cdot zx\phi^{-1}$. Furthermore, it is shown that for all $x\in L$, $H(L)$ is an Osborn loop if and only if $(L,\cdot)$ is an Osborn loop, $(x\alpha\cdot x^{\rho})x=x\alpha$, $x(x^{\lambda}\cdot x\phi^{-1})=x\phi^{-1}$ and every pair of automorphisms in $A(L)$ is nuclear (i.e. $x\alpha\cdot x^{\rho},x^{\lambda}\cdot x\phi\in N(L,\cdot )$). It is shown that if $H(L)$ is an Osborn loop, then $A(L,\cdot)= \mathcal{P}(L,\cdot)\cap\Lambda(L,\cdot)\cap\Phi(L,\cdot)\cap\Psi(L,\cdot)$
and for any $\alpha\in A(L)$, $\alpha= L_{e\pi}=R^{-1}_{e\varrho}$ for some
$\pi\in \Phi(L,\cdot)$ and some $\varrho\in \Psi(L,\cdot)$. Some commutative diagrams are deduced by considering isomorphisms among the various groups of regular bijections (whose intersection is $A(L)$) and the nucleus of $(L,\cdot)$.
\end{abstract}
\section{Introduction}
\paragraph{}
The holomorph of a loop is a loop according to Bruck \cite{phd15}. Since then, the concept of holomorphy of loops
has caught the attention of to some researchers. Interestingly, Adeniran $\cite{phd11}$ and Robinson \cite{phd40},
Chein and Robinson\cite{phd17}, Adeniran et. al. \cite{phd36}, Chiboka and Solarin $\cite{phd19}$,
$\cite{phd20}$, Bruck \cite{phd15}, Bruck and Paige \cite{phd16}, Robinson
\cite{phd39}, Huthnance \cite{phd24} have
respectively studied the holomorphic structures of Bol/Bruck loops, Moufang loops, central loops,
conjugacy closed loops, inverse property loops, A-loops, extra
loops and weak inverse property loops.

After the discovery of Osborn loops by Osborn \cite{phd35} and Huthnance \cite{phd24}, Osborn loops were formally introduced and studied by Basarab \cite{phd148,phd46,phd14,phd137,phd13}, in the $\textrm{20}^{\textrm{th}}$ century. In this $\textrm{21}^{\textrm{st}}$ century, the study of Osborn loops was recently revived by Kinyon \cite{phd29}, where he proposed some problems and heart burning questions. Some of this problems and questions have been solved, answered fully or partially or somewhat addressed in Jaiy\'e\d ol\'a et. al. \cite{phd195,phd195d,open4}, Jaiy\'e\d ol\'a and Ad\'en\'iran \cite{phd27,phd28,phd195c} and Jaiy\'e\d ol\'a \cite{phd195b,phd26,open6,osbornnew}. Some results on the application of Osborn loops to cryptography can be found in Jaiy\'e\d ol\'a and Ad\'en\'iran \cite{phd195e} and Jaiy\'e\d ol\'a \cite{phd195f,phd195g}.

Some popular varieties of Osborn loops are: extra loops, Moufang loops, CC-loops, universal
WIPLs and V.D. loops. Some studies on them can be found in Dr\'apal \cite{phd21,phd22,phd23,phd38,phd139}, Cs\"org\H o and Dr\'apal \cite{phd106}, Cs\"org\H o \cite{phd108}, Kinyon and Kunen \cite{phd36.1,phd47}, Kinyon et. al. \cite{phd30}. Some newly constructed Osborn loops can be found in Isere et. al. \cite{phd25,phd25b,open5.1,open5.2,open5.4}, Adeniran and Isere \cite{isereref}.
\section{Preliminaries}
\paragraph{}
Let $G$ be a non-empty set. Define a binary operation ($\cdot $) on
$G$. If $x\cdot y\in G$ for all $x, y\in G$, then the pair $(G, \cdot )$
is called a groupoid or Magma.

If each of the equations:
\begin{displaymath}
a\cdot x=b\qquad\textrm{and}\qquad y\cdot a=b
\end{displaymath}
has unique solution in $G$ for $x$ and $y$ respectively, then $(G,
\cdot )$ is called a quasigroup.

If there exists a unique element $e\in G$ called the
identity element such that for all $x\in G$, $x\cdot
e=e\cdot x=x$, $(G, \cdot )$ is called a loop. We write
$xy$ instead of $x\cdot y$, and stipulate that $\cdot$ has lower
priority than juxtaposition among factors to be multiplied. For
instance, $x\cdot yz$ stands for $x(yz)$.

For a groupoid $(G, \cdot )$, the right translation of $x$ i.e. $R_x~:G\to G$ is defined by
$yR_x=y\cdot x$ while the left translation of $x$ i.e. $L_x~:G\to G$ is
defined by $yL_x=x\cdot y$ for all $x,y\in G$.

It can now be seen that a groupoid $(G, \cdot )$ is a quasigroup if
its left and right translation mappings are bijections or
permutations. Since the left and right translation mappings of a
loop are bijective, then the inverse mappings $L_x^{-1}$ and
$R_x^{-1}$ exist. Let
\begin{displaymath}
x\backslash y =yL_x^{-1}\qquad\textrm{and}\qquad
x/y=xR_y^{-1}
\end{displaymath}
and note that
\begin{displaymath}
x\backslash y =z\Longleftrightarrow x\cdot
z=y\qquad\textrm{and}\qquad x/y=z\Longleftrightarrow z\cdot y=x.
\end{displaymath}
Hence, $(G, \backslash )$ and $(G, /)$ are also quasigroups. Using
the operations ($\backslash$) and ($/$), the definition of a loop
can be stated as follows.

\begin{mydef}\label{0:1}
A \textit{loop} $(G,\cdot ,/,\backslash ,e)$ is a set $G$ together
with three binary operations ($\cdot $), ($/$), ($\backslash$) and
one nullary operation $e$ such that
\begin{description}
\item[(i)] $x\cdot (x\backslash y)=y$, $(y/x)\cdot x=y$ for all
$x,y\in G$,
\item[(ii)] $x\backslash (x\cdot y)=y$, $(y\cdot x)/x=y$ for all
$x,y\in G$ and
\item[(iii)] $x\backslash x=y/y$ or $e\cdot x=x$ for all
$x,y\in G$.
\end{description}
\end{mydef}
We also stipulate that ($/$) and ($\backslash$) have higher priority
than ($\cdot $) among factors to be multiplied. For instance,
$x\cdot y/z$ and $x\cdot y\backslash z$ stand for $x(y/z)$ and
$x\cdot (y\backslash z)$ respectively.

In a loop $(G,\cdot )$ with identity element $e$, the left
inverse element of $x\in G$ is the element $xJ_\lambda
=x^\lambda\in G$ such that
\begin{displaymath}
x^\lambda\cdot x=e
\end{displaymath}
while the right inverse element of $x\in G$ is the element
$xJ_\rho =x^\rho\in G$ such that
\begin{displaymath}
x\cdot x^\rho=e.
\end{displaymath}

A loop is called an Osborn loop if it obeys any of the three identities
\begin{equation}
x(yz\cdot x)= (x^{\lambda}\backslash y)\cdot zx
\end{equation}
\begin{equation}
x(yz\cdot x)= x(yx^{\lambda}\cdot x)\cdot zx
\end{equation}
\begin{equation}
x(yz\cdot x)= x(yx\cdot x^{\rho})\cdot zx
\end{equation}
Given any two sets $X$ and $Y$. The statement '$f~:X\rightarrow Y$ is defined as $f(x)=y,~x\in X,~y\in Y$' will be expressed as '$f~:X\rightarrow Y~\uparrow~f(x)=y$'.

Let $(G,\cdot )$ be a loop and let $A,B$ and $C$ be three bijective
mappings, that map $G$ onto $G$. The identity mapping on $G$ will be denoted by $I$. The triple $\alpha =(A,B,C)$ is
called an autotopism of $(G,\cdot )$ if and only if
\begin{displaymath}
xA\cdot yB=(x\cdot y)C~\forall~x,y\in G.
\end{displaymath}
Such triples form a group
$AUT(G,\cdot )$ called the autotopism group of $(G,\cdot )$ under the binary operation of componentwise composition. That is, for $(A_1,B_1,C_1),(A_2,B_2,C_2)\in AUT(G,\cdot ),~(A_1,B_1,C_1)(A_2,B_2,C_2)=(A_1A_2,B_1B_2,C_1C_2)$.

If $A=B=C$, then $A$ is called an automorphism of the
loop $(G,\cdot )$. Such bijections form a
group $AUM(G,\cdot )$ called the automorphism group of $(G,\cdot )$. Let $G$ and $H$ be groups such that $\varphi:G\to H$ is an isomorphism. If $\varphi (g)=h$, then this would be expressed as $g\overset{\varphi}{\cong}h$.

\begin{mydef}
Let $(Q,\cdot)$ be a loop and $A(Q)\le AUM(Q,\cdot)$
be a group of automorphisms of the loop $(Q,\cdot)$. Let $H=A(Q)\times Q$. Define $\circ$ on $H$ as
\begin{displaymath}
(\alpha,x)\circ(\beta,y)=(\alpha\beta,x\beta\cdot y)~\textrm{for all}~(\alpha,x),(\beta,y)\in H.
\end{displaymath}
$(H,\circ)$ is a loop and is called the A-holomorph of $(Q,\cdot)$.
\end{mydef}

The right nucleus of $(L,\cdot)$ is defined by $N_\rho (L,\cdot)=\{x\in L~|~zy\cdot x=z\cdot yx~\forall~y,z\in L\}$.
The left nucleus of $(L,\cdot)$ is defined by $N_\lambda (L,\cdot)=\{x\in L~|~x\cdot yz=xy\cdot z~\forall~y,z\in L\}$.
The middle nucleus of $(L,\cdot)$ is defined by $N_\mu (L,\cdot)=\{x\in L~|~zx\cdot y=z\cdot xy~\forall~y,z\in L\}$.
The nucleus of $(L,\cdot)$ is defined by $N(L,\cdot)=N_\rho (L,\cdot)\cap N_\lambda (L,\cdot)\cap N_\mu (L,\cdot)$.
The centrum of $(L,\cdot)$ is defined by $C(L, \cdot )=\{a\in L : ax=xa~\forall~x\in L\}$ while its center is defined by
$Z(L, \cdot )=N(L, \cdot )\cap C(L, \cdot )$.

Let $(G,\cdot )$ be a quasigroup. Then
\begin{enumerate}
\item a bijection $U$ is called autotopic if there exists $(U,V,W)\in AUT(G,\cdot )$; the set of all such mappings forms a group $\Sigma(G,\cdot )$.
\item a bijection $U$ is called $\rho$-regular if there exists $(I,U,U)\in AUT(G,\cdot )$; the set of all such mappings forms a group $\mathcal{P}(G,\cdot )$.
    \item a bijection $U$ is called $\lambda$-regular if there exists $(U,I,U)\in AUT(G,\cdot )$; the set of all such mappings forms a group $\Lambda(G,\cdot )\le\Sigma(G,\cdot )$.
\item a bijection $U$ is called $\mu$-regular if there exists a bijection $U'$ such that $(U,U'^{-1},I)\in AUT(G,\cdot )$. $U'$ is called the adjoint of $U$. The set of all $\mu$-regular mappings forms a group $\Phi(G,\cdot )\le\Sigma(G,\cdot )$. The set of all adjoint mapping
    forms a group $\Psi(G,\cdot )$.
\end{enumerate}
\begin{myth}(Jaiy\'e\d ol\'a  \cite{phd41})
Let $(G,\cdot )$ be a loop. Let
\begin{displaymath}
\psi : \mathcal{P}(G,\cdot)\rightarrow N_{\rho}(G,\cdot)\uparrow \psi(U)=eU,\delta : \Lambda(G,\cdot)\rightarrow N_{\lambda}(G,\cdot)\uparrow \delta(U)=eU,~\varphi :\Phi(G,\cdot)\rightarrow \Psi(G,\cdot)
\end{displaymath}
\begin{displaymath}
\uparrow \varphi(U)=U',~\sigma :\Phi(G,\cdot)\rightarrow N_{\mu}(G,\cdot)\uparrow \sigma(U)=eU~ \textrm{and}~ \beta: \Psi(G,\cdot)\rightarrow N_{\mu}(G,\cdot)\uparrow \beta(U')=eU'
\end{displaymath}
Then $\mathcal{P}(G,\cdot )\stackrel{\psi}{\cong} N_{\rho}(G,\cdot ),~\Lambda(G,\cdot )\stackrel{\delta}{\cong} N_{\lambda}(G,\cdot ),~\Phi(G,\cdot )\stackrel{\varphi}{\cong}\Psi(G,\cdot ),~\Phi(G,\cdot )\stackrel{\sigma}{\cong}N_\mu(G,\cdot ),~\Psi(G,\cdot )\stackrel{\beta}{\cong}N_\mu(G,\cdot )$.
\end{myth}

\section{Main Results}
\subsection{Holomorph of an Osborn loop}
\begin{myth}\label{isere1}
Let $(L,\cdot)$ be a loop and A(L) be a group of automorphisms of $(L,\cdot)$. Then, the
A(L)-holomorph $(H,\circ)$ of $(L,\cdot)$ is an Osborn loop if and only if
\begin{equation}
x\alpha (yz\cdot x\phi^{-1})= x\alpha (yx^\lambda\cdot x) \cdot zx\phi^{-1} ~\forall~x,y,z\in L ~\textrm{and}~ \alpha, \phi\in A(L)
\end{equation}
\end{myth}
{\bf Proof:}\\
Suppose A(L)-holomorph $(H,\circ)$ of $(L,\cdot)$ is an Osborn loop, then we have
\begin{displaymath}
(\alpha, x)\circ \{[(\beta, y)\circ (\gamma, z)]\circ (\alpha, x)\}= (\alpha, x)\circ \{[(\beta,y)\circ (\alpha, x)^{\lambda}]\circ (\alpha, x)\}\{(\gamma, z)\circ(\alpha, x)\}
\end{displaymath}
\begin{displaymath}
\Leftrightarrow\big(\alpha(\beta\gamma\alpha), x\beta\gamma\alpha [(y\gamma \cdot z)\alpha\cdot x]\big)= (\alpha, x)\circ \{(\beta\alpha^{-1}, yx^\lambda \cdot x)\circ (\alpha, x)\}\circ (\gamma\alpha, z\alpha \cdot x)
\end{displaymath}
\begin{displaymath}
\Leftrightarrow\big(\alpha(\beta\gamma\alpha), x\beta\gamma\alpha [(y\gamma \cdot z)\alpha\cdot x]\big)= [(\alpha, x)\circ (\beta, yx^\lambda\cdot x)]\circ (\gamma\alpha, z\alpha \cdot x)
\end{displaymath}
\begin{displaymath}
\Leftrightarrow \big(\alpha(\beta\gamma\alpha), x\beta\gamma\alpha [(y\gamma \cdot z)\alpha\cdot x]\big)= \big(\alpha(\beta\gamma\alpha), (x\beta (yx^\lambda\cdot x))\gamma\alpha\cdot (z\alpha\cdot x)\big)
\end{displaymath}
if and only if
\begin{displaymath}
 x\beta\gamma\alpha[(y\gamma \cdot z)\alpha\cdot x]=  (x\beta \cdot (yx^\lambda\cdot x))\gamma\alpha\cdot  (z\alpha \cdot x) \Leftrightarrow
\end{displaymath}
\begin{displaymath}
 x\beta\gamma\alpha \big[(y\gamma\alpha \cdot z\alpha) \cdot x\big]=  \big(x\beta\gamma\alpha ((yx^\lambda)\gamma\alpha\cdot x\gamma\alpha)\big)(z\alpha \cdot x) \qquad \forall~ x,y,z\in L ~\textrm{and}~\alpha,\beta, \gamma \in A(L).
\end{displaymath}
Putting $\phi = \gamma\alpha$ , we have
\begin{displaymath}
 x\beta\phi[(y\phi \cdot z\alpha)\cdot x]=  \big(x\beta\phi \cdot ((yx^\lambda)\phi \cdot x\phi)\big)(z\alpha \cdot x) \qquad \forall~ x,y,z\in L ~\textrm{and}~\alpha,\beta, \phi \in A(L).
\end{displaymath}
Therefore,
\begin{displaymath}
 x\beta\big[(y \cdot z\alpha\phi^{-1})x\phi^{-1}\big]=  (x\beta \cdot (yx^\lambda \cdot x))(z\alpha\phi^{-1} \cdot x\phi^{-1}).
\end{displaymath}
Letting $\bar{x}=x\phi^{-1}$ and $x = \bar{x}\phi$, $\bar{z}=z\alpha\phi^{-1}$, we obtain
\begin{displaymath}
 \bar{x}\phi\beta(y\bar{z} \cdot \bar{x})=  \bar{x}\phi\beta(y{(\bar{x}\phi)}^\lambda \cdot \bar{x}\phi) \cdot  \bar{z}\bar{x}
\end{displaymath}
Again, since $\phi$ is an automorphism, then letting $\bar{x}\phi = x$ and $\bar{x}=x\phi^{-1}$, and replacing $\bar{z}$ with $z$ and $\beta$ with $\alpha$, we obtain
\begin{displaymath}
 x\alpha(yz \cdot x\phi^{-1})=  x\alpha(yx^\lambda \cdot x) \cdot  zx\phi^{-1} ~ \forall~ x,y,z\in L ~\textrm{and}~\alpha,\phi\in A(L).
\end{displaymath}
The converse is obtained by reversing the process.

\begin{mycor}\label{isere2}
Let $(L,\cdot)$ be a loop, and A(L) be the group of all automorphisms of $L$. Then, the
holomorph $(H,\circ)$ of $(L,\cdot)$ is an Osborn loop if and only if
\begin{displaymath}
(R_{x^\lambda}R_{x}L_{x\alpha}, R_{x\phi^{-1}}, R_{x\phi^{-1}}L_{x\alpha})
\end{displaymath}
is an autotopism of $L$ for all $x\in L$ and all $\alpha, \phi \in A(L)$.
\end{mycor}
{\bf Proof:}\\
This is a consequence of Theorem 3.1.
\begin{mylem}\label{isere3}
Let A(L) be an automorphism group of an Osborn loop $(L,\cdot)$. The holomorph
$(H,\circ)$ of $(L,\cdot)$ is Osborn if and only if the triples
\begin{displaymath}
(L_{x}^{-1}L_{x\alpha}, I, L_{x}^{-1}L_{x\alpha}) ~\textrm{and}~
(I,R_{x}^{-1}R_{x\phi^{-1}},L_{x}^{-1}R_{x}^{-1}R_{x\phi^{-1}}L_{x})
\end{displaymath}
are autotopisms of $L$ for all $x \in L$ and all $\alpha \in A(L)$.
\end{mylem}
{\bf Proof:}\\
\begin{equation}\label{isere18}
\textrm{Let}~ A= (R_{x^\lambda}R_{x}L_{x\alpha}, R_{x\phi^{-1}}, R_{x\phi^{-1}}L_{x\alpha})~\textrm{and}~B= (R_{x^\lambda}R_{x}L_{x}, R_{x}, R_{x}L_{x})
\end{equation}
Since $(L,\cdot)$ is an Osborn loop, $B$ is an autotopism of $L$ for all $x\in L$. The holomorph
of $(L,\cdot)$ is an Osborn if and only if A is autotopism of $L$ (by Corollary~\ref{isere2}). So, the triple
\begin{equation}\label{isere20}
 B^{-1}=(L_{x}^{-1}R_{x}^{-1}R_{x^\lambda}^{-1}, R_{x}^{-1}, L_{x}^{-1}R_{x}^{-1})
\end{equation}
is also an autotopism of $L$ for all $x \in L$. Hence, $(H,\circ)$ is an Osborn loop if and only if
\begin{equation}\label{isere21}
B^{-1}A=(L_{x}^{-1}L_{x\alpha}, R_{x}^{-1}R_{x\phi^{-1}}, L_{x}^{-1}R_{x}^{-1}R_{x\phi^{-1}}L_{x\alpha})\in AUT(L,\cdot ).
\end{equation}
\begin{displaymath}
\textrm{Thus,}~yL_{x}^{-1}L_{x\alpha}\cdot zR_{x}^{-1}R_{x\phi^{-1}} =(yz)L_{x}^{-1}R_{x}^{-1}R_{x\phi^{-1}}L_{x\alpha}
\end{displaymath}
\begin{equation}\label{isere19}
\Leftrightarrow[(x\alpha)\cdot (x\backslash y)]\cdot [(z/x)\cdot x\phi^{-1}]= (x\alpha)\cdot
\{[x\backslash (yz)]/x \cdot x\phi^{-1}\}.
\end{equation}
Put $\phi = I$ into equation~\eqref{isere19} to get
\begin{equation}\label{isere17}
[(x\alpha)\cdot (x\backslash y)]\cdot z = [(x\alpha)\cdot
(x\backslash (yz)]
\end{equation}
\begin{equation*}\label{isere15}
\Leftrightarrow y L_{x}^{-1}L_{x\alpha}\cdot z =(yz)L_{x}^{-1}L_{x\alpha}\Leftrightarrow(L_{x}^{-1}L_{x\alpha},I,L_{x}^{-1}L_{x\alpha})\in AUT(L,\cdot ).
\end{equation*}
Now, putting $\alpha = I$ into equation~\eqref{isere19}, we obtain
\begin{gather*}\label{isere13}
[x\cdot (x\backslash y)]\cdot [(z/x)\cdot x\phi^{-1}]= x \cdot \{[x\backslash(yz)]/x \cdot
x\phi^{-1}\}\\
\Leftrightarrow y \cdot zR_{x}^{-1}R_{x\phi^{-1}}= \{[x\backslash(yz)]/x \cdot
x\phi^{-1}\}L_{x}
\Leftrightarrow(I,R_{x}^{-1}R_{x\phi^{-1}},L_{x}^{-1}R_{x}^{-1}R_{x\phi^{-1}}L_{x})\in AUT(L,\cdot ).
\end{gather*}
The converse follows from Theorem ~\ref{isere1} and Corollary~\ref{isere2}.
\begin{myth}\label{isere4}
Let A(L) be an automorphism group of a loop $(L,\cdot)$. The holomorph $(H,\circ)$  of  $(L,\cdot)$ is an Osborn loop if and only if:
\begin{itemize}
\item[(i)]
 $(L,\cdot)$ is an Osborn loop,
\item[(ii)]
$x\alpha \cdot x^{\rho}, x^{\lambda}\cdot x\phi \in N(L,\cdot)$,
\item[(iii)]
$(x\alpha\cdot x^{\rho})x = x\alpha$,
\item[(iv)]
$x(x^{\lambda} \cdot x\phi^{-1})=x\phi^{-1}$,
\end{itemize}
for every $x,y \in L$ and $\alpha, \phi \in A(L)$.
 \end{myth}
{\bf Proof:}\\
\begin{itemize}
\item[(i)] Suppose $(H,\circ)$ is an Osborn loop. $(K,\circ)$ is a subloop of $(H,\circ)$ given by
$K=\{(I,x): x\in L\}$.
Therefore, $(L,\cdot)\cong(K,\circ)$. Since $(K,\circ)$ is
an Osborn loop, it follows that $(L,\cdot)$ is an Osborn loop.
\item[(ii)] Since $(H,\circ)$ is Osborn, then by Lemma~\ref{isere3}, we have
$(L_{x}^{-1}L_{x\alpha},I,L_{x}^{-1}L_{x\alpha})\in AUT(L,\cdot )$
\begin{equation}\label{isere24}
\Leftrightarrow yL_{x}^{-1}L_{x\alpha} \cdot  zI = (yz)L_{x}^{-1}L_{x\alpha}~\textrm{for all}~y,z \in (L,\cdot).
\end{equation}
Putting $y = e$ in equation~\eqref{isere24} gives:
$eL_{x}^{-1}L_{x\alpha} \cdot  zI = (ez)L_{x}^{-1}L_{x\alpha}\Rightarrow
(x\alpha\cdot x^{\rho})z = (x\alpha)(x\backslash z)$
\begin{equation}\label{isere24.11}
\Rightarrow
L_{x\alpha\cdot x^{\rho}}  = L_{x}^{-1}L_{x\alpha}.
\end{equation}
Again, since $(I,R_{x}^{-1}R_{x\phi^{-1}},L_{x}^{-1}R_{x}^{-1}R_{x\phi^{-1}}L_{x})\in AUT(L,\cdot )$, we have
\begin{equation}\label{isere23}
yI\cdot zR_{x}^{-1}R_{x\phi^{-1}} = (yz)L_{x}^{-1}R_{x}^{-1}R_{x\phi^{-1}}L_{x},
\end{equation}
putting $z = e$, we have
\begin{displaymath}
y\cdot eR_{x}^{-1}R_{x\phi^{-1}} = (ye)L_{x}^{-1}R_{x}^{-1}R_{x\phi^{-1}}L_{x}\Rightarrow y\cdot (x^{\lambda}\cdot x\phi^{-1}) = yL_{x}^{-1}R_{x}^{-1}R_{x\phi^{-1}}L_{x}
\end{displaymath}
\begin{equation}
\Rightarrow R_{x^{\lambda }\cdot x\phi^{-1}} = L_{x}^{-1}R_{x}^{-1}R_{x\phi^{-1}}L_{x}.
\end{equation}
Also, substituting $y = e$ in equation~\eqref{isere23}, we get
\begin{equation}\label{isere25.1}
R^{-1}_{x}R_{x\phi^{-1}} = L_{x}^{-1}R_{x}^{-1}R_{x\phi^{-1}}L_{x}
\end{equation}
\begin{equation}\label{isere25}
\textrm{so},~R_{x^{\lambda }\cdot x\phi^{-1}} = R^{-1}_{x}R_{x\phi^{-1}} =L_{x}^{-1}R_{x}^{-1}R_{x\phi^{-1}}L_{x}.
\end{equation}
So, $x\alpha\cdot x^{\rho}\in N_{\lambda}(L,\cdot)$ and $x^{\lambda} \cdot x\phi\in N_{\rho}(L,\cdot)$, hence, $x^{\lambda} \cdot x\phi,x\alpha\cdot x^{\rho}\in N(L,\cdot)$ for
all $x \in L$ and all $\alpha,\phi \in A(L)$.
\item[(iii)] From equation~\eqref{isere24.11},
\begin{displaymath}
L_{x}L_{x\alpha \cdot x^{\rho}}= L_{x\alpha}
\Leftrightarrow
(x\alpha\cdot x^{\rho})\cdot xy = x\alpha \cdot y.
\end{displaymath}
Since $(x\alpha\cdot x^{\rho})\in N_{\lambda}(L,\cdot)$, then
\begin{equation}
(x\alpha\cdot x^{\rho})x=x\alpha~\textrm{for
all}~x \in L,~\alpha,\phi \in A(L).
\end{equation}
\item[(iv)] From equation~\eqref{isere25},
\begin{equation*}
R_xL_xR_{x^{\lambda}\cdot x\phi^{-1}} = R_{x\phi^{-1}}L_x\Leftrightarrow
(x\cdot yx)(x^{\lambda}\cdot x\phi^{-1})=x(y\cdot x\phi^{-1}).
\end{equation*}
Since $x^{\lambda}\cdot x\phi^{-1}\in N_{\rho}(L,\cdot)$, $x\cdot (yx)(x^{\lambda}\cdot x\phi^{-1})=x(y\cdot x\phi^{-1})\Rightarrow$
\begin{equation}
   x(x^{\lambda}\cdot x\phi^{-1})= x\phi^{-1}.
\end{equation}
\end{itemize}
The converse:
suppose $(L,\cdot)$ is an Osborn loop such that (ii), (iii) and (iv) hold. We need
to show that $(H,\circ)$ is an Osborn loop.\\
Already, $(x\alpha\cdot x^{\rho})x=x\alpha$, thence $(x\alpha\cdot x^{\rho})x\cdot y=x\alpha\cdot y$.

Since $x\alpha\cdot x^{\rho} \in N_{\lambda}(L,\cdot)$, $(x\alpha\cdot x^{\rho})\cdot xy=(x\alpha\cdot x^{\rho})x\cdot y=x\alpha\cdot y\Rightarrow
L_{x\alpha\cdot x^{\rho}}=L^{-1}_{x}L_{x\alpha}$.
Next, since $(x^{\lambda}\cdot x\phi^{-1})\in N_{\rho}(L,\cdot)$, then
\begin{gather*}
x(x^{\lambda}\cdot x\phi^{-1})=x\phi^{-1}\Rightarrow (yx)(x^{\lambda}\cdot x\phi^{-1})=(y\cdot x\phi^{-1})\\
\Rightarrow (x\cdot yx)(x^{\lambda}\cdot x\phi^{-1})=x(y\cdot x\phi^{-1})\Rightarrow R_{x^\lambda}\cdot x\phi^{-1}=L^{-1}_{x}R^{-1}_{x}R_{x\phi^{-1}}L_{x}.
\end{gather*}
Already, $x(x^\lambda\cdot x\phi^{-1})=x\phi^{-1}$. Since, $x^\lambda\cdot x\phi^{-1}\in N(L,\cdot )$, then $y\cdot x(x^\lambda\cdot x\phi^{-1})=y\cdot x\phi^{-1}\Rightarrow yx\cdot (x^\lambda\cdot x\phi^{-1})=y\cdot x\phi^{-1}\Rightarrow $
\begin{equation}
R_{x^\lambda\cdot x\phi^{-1}}=R^{-1}_{x}R_{x\phi^{-1}}
\end{equation}
Since $x\alpha\cdot x^{\rho}\in N_{\lambda}(L,\cdot)$, then $(L^{-1}_xL_{x\alpha},I,L^{-1}_xL_{x\alpha})\in AUT(L,\cdot )$.

And also, since $x^\lambda\cdot x\phi^{-1}\in N_\rho(L,\cdot)$, then $(I,R^{-1}_xR_{x\phi^{-1}},L^{-1}_{x}R^{-1}_{x}R_{x\phi^{-1}}L_{x})\in AUT(L,\cdot )$. Hence, by Lemma~\ref{isere3}, the holomorph $(H,\circ)$ of $(L,\cdot)$
is an Osborn loop.
\begin{mylem}\label{isere6}
Let $A(L)$ be an automorphism group of a loop $(L,\cdot)$. If the holomorph $(H,\circ)$
of $(L,\cdot)$ is an Osborn loop, then the following identities hold:
\begin{itemize}
\item[(1)] $(x\alpha\cdot x^{\rho})\cdot xy= x\alpha\cdot y$; $x\cdot (x\alpha)^\rho=(x\alpha\cdot x^{\rho})^{\rho}$,~ $(x\alpha\cdot x^\rho)x=x\alpha$,
\item[(2)] $(x\cdot yx)(x^\lambda \cdot x\phi^{-1})= x(y\cdot x\phi^{-1})$; $(x\phi^{-1})^\lambda \cdot x=(x^\lambda \cdot x\phi^{-1})^\lambda$,
\item[(3)] $yx\cdot (x^\lambda \cdot x\phi^{-1})=y\cdot x\phi^{-1}$; $x(x^\lambda \cdot x\phi^{-1})=x\phi^{-1}$,
    \item[(4)] $x(y/x^\lambda \cdot x\phi^{-1})=(xy)/x\cdot x\phi^{-1}$; $x(x^\rho/x^\lambda \cdot x\phi^{-1})=x^\lambda\cdot x\phi^{-1}$,
    \end{itemize}
    for all $x,y \in L$ and $\alpha,\phi \in A(L)$.
\end{mylem}
{\bf Proof:}\\
From Theorem~\ref{isere4}, we have
\begin{itemize}
\item[(1)]
\begin{equation}\label{ise1}
L_{x\alpha\cdot x^\rho}=L^{-1}_{x}L_{x\alpha}\Rightarrow(x\alpha\cdot x^\rho)\cdot xy = x\alpha\cdot y
\end{equation}
Put $y=(x\alpha)^\rho$ in \eqref{ise1} to get $(x\alpha\cdot x^{\rho})^{\rho}=x\cdot (x^\alpha)^\rho$. Putting $y=e$ in \eqref{ise1}, then $(x\alpha\cdot x^\rho)x=x\alpha$.
\item[(2)]
\begin{equation}\label{ise2}
R_{x^{\lambda}\cdot x\phi^{-1}}= L^{-1}_xR^{-1}_xR_{x\phi^{-1}}L_x\Rightarrow
(x\cdot yx)(x^{\lambda}\cdot x\phi^{-1})=x(y\cdot x\phi^{-1})
\end{equation}
Put $y=(x\phi^{-1})^\lambda$ in \eqref{ise2} to get $(x^{\lambda}\cdot x\phi^{-1})^\lambda =(x\phi^{-1})^\lambda x$.
\item[(3)]
\begin{equation}\label{ise3}
R_{x^{\lambda}\cdot x\phi^{-1}}= R^{-1}_xR_{x\phi^{-1}}
\Rightarrow yx\cdot(x^{\lambda}\cdot x\phi^{-1})=y\cdot x\phi^{-1}
\end{equation}
Put $y=e$ in \eqref{ise2} to get $x(x^\lambda \cdot x\phi^{-1})=x\phi^{-1}$.
\item[(4)] $R^{-1}_{x}R_{x\phi^{-1}} =L_{x}^{-1}R_{x}^{-1}R_{x\phi^{-1}}L_{x}\Rightarrow L_{x}R^{-1}_{x}R_{x\phi^{-1}} =R_{x}^{-1}R_{x\phi^{-1}}L_{x}\Rightarrow x(y/x^\lambda \cdot x\phi^{-1})=(xy)/x\cdot x\phi^{-1}$. Put $y=x^\rho$, then $x(x^\rho/x^\lambda \cdot x\phi^{-1})=x^\lambda\cdot x\phi^{-1}$.
\end{itemize}
The proof is complete.
\begin{mylem}\label{isere7}
Let $(L,\cdot)$ be a loop. If the A-holomorph $H(L)$ of $L$ is an Osborn loop, then for all $x\in L$ and $\alpha,\phi \in A(L)$.
\begin{itemize}
\item[(a)]
$(L^{-1}_xL_{x\alpha},I,L^{-1}_xL_{x\alpha}) \in AUT(L,\cdot)$.
\item[(b)]
$(L_{x\alpha\cdot x^\rho},I,L_{x\alpha\cdot x^\rho})\in AUT(L,\cdot)$.
\item[(c)]
$(I,R^{-1}_xR_{x\phi^{-1}},R^{-1}_xR_{x\phi^{-1}})\in AUT(L,\cdot)$.
\item[(d)]
$(I,L^{-1}_xR^{-1}_xR_{x\phi^{-1}}L_x,L^{-1}_xR^{-1}_xR_{x\phi^{-1}}L_x)\in AUT(L,\cdot)$.
\item[(e)]
$(I,R_{x^{\lambda}\cdot x\phi^{-1}},R_{x^{\lambda}\cdot x\phi^{-1}})\in AUT(L,\cdot)$.
\item[(f)]
$(R_{x\alpha\cdot x^\rho},L^{-1}_{x\alpha\cdot x^\rho},I), (R_{x^{\lambda}\cdot x\phi^{-1}},L^{-1}_{x^{\lambda}\cdot x\phi^{-1}},I) \in AUT(L,\cdot)$.
\end{itemize}
\end{mylem}
{\bf Proof:}\\
\begin{itemize}
\item[(a)] Following the steps in the proof of Lemma~\ref{isere3}, we obtain (a).
\item[(b)] Use Theorem~\ref{isere4} and the fact that $L_{x\alpha\cdot x^\rho}=L^{-1}_xL_{x\alpha}$.
\item[(c)] Follow the steps in Lemma~\ref{isere3} and Theorem~\ref{isere4}.
\item[(d)] Follow the steps in Lemma~\ref{isere3} and Theorem~\ref{isere4}.
\item[(e)] Follow the steps in Lemma~\ref{isere3} and Theorem~\ref{isere4}.
\item[(f)] Since $x\alpha\cdot x^{\rho}\in N(L,\cdot)$, obviously, it is in $N_{\mu}(L,\cdot)$.
Then:
\begin{equation}
x(x\alpha\cdot x^\rho)\cdot y=x\cdot (x\alpha\cdot x^\rho)y\Rightarrow
xR_{x\alpha\cdot x^\rho}\cdot yL^{-1}_{x\alpha\cdot x^\rho} = xy
\end{equation}
which implies that
\begin{displaymath}
(R_{x\alpha\cdot x^\rho},L^{-1}_{x\alpha\cdot x^\rho}, I) \in AUT(L,\cdot).
\end{displaymath}
Since $x^{\lambda}\cdot x\phi^{-1} \in N(L,\cdot)\Rightarrow x^{\lambda}\cdot x\phi^{-1} \in N_{\mu}(L,\cdot)$, then by definition,
\begin{equation*}
x(x^{\lambda}\cdot x\phi^{-1})\cdot y =x\cdot (x^{\lambda}\cdot x\phi^{-1})y\Rightarrow
xR_{x^{\lambda}\cdot x\phi^{-1}}\cdot y =x\cdot yL_{x^{\lambda}\cdot x\phi^{-1}}\Rightarrow
\end{equation*}
\begin{displaymath}
(R_{x^{\lambda}\cdot x\phi^{-1}},L_{x^{\lambda}\cdot x\phi^{-1}},I) \in AUT(L,\cdot).
\end{displaymath}
\end{itemize}
That completes the proof.
\begin{mycor}\label{isere3.3}
Let $(L,\cdot)$ be a loop. If the A-holomorph $H(L)$ of $L$ is an Osborn loop, then for all $x\in L$ and $\alpha,\phi \in A(L)$.
\begin{itemize}
\item[(a)]
$L^{-1}_xL_{x\alpha}, L_{x\alpha\cdot x^\rho}\in \Lambda(L,\cdot); L_{x\alpha}\in L_x\Lambda(L,\cdot)$.
\item[(b)]
$R^{-1}_xR_{x\phi^{-1}},L^{-1}_xR^{-1}_xR_{x\phi^{-1}}L_x,
R_{x^{\lambda}\cdot x\phi^{-1}}\in \mathcal{P}(L,\cdot)$; $R_{x\phi^{-1}}\in R_x\mathcal{P}(L,\cdot), R_{x\phi^{-1}}L_x\in R_xL_x\mathcal{P}(L,\cdot)$.
\item[(c)]
$R_{x\alpha\cdot x^\rho}, R_{x^{\lambda}\cdot x\phi^{-1}}\in \Phi(L,\cdot)$, $L_{x\alpha\cdot x^\rho}, L_{x^{\lambda}\cdot x\phi^{-1}}\in \Psi(L,\cdot)$.
\end{itemize}
\end{mycor}
{\bf Proof:}\\
Use Lemma~\ref{isere7}.
\begin{myth}\label{isere3.4}
Let $L$ be a loop and $H(L)$ its A-holomorph. If $H(L)$ is an Osborn loop, then
$A(L,\cdot)= \mathcal{P}(L,\cdot)\cap\Lambda(L,\cdot)\cap\Phi(L,\cdot)\cap\Psi(L,\cdot)$
and for any $\alpha\in A(L)$, $\alpha= L_{e\pi}=R^{-1}_{e\varrho}$ for some
$\pi\in \Phi(L,\cdot)$ and some $\varrho \in \Psi(L,\cdot)$
\end{myth}
{\bf Proof:}\\
From Corollary~\ref{isere3.3}, $L_{x\alpha}\in L_x\Lambda(L,\cdot)\Rightarrow L_{x\alpha}=L_x\lambda$ for some $\lambda\in \Lambda(L,\cdot)$. So
\begin{equation}
x\alpha\cdot y= (xy)\lambda
\end{equation}
implies that,
$( \alpha,I,\lambda)\in AUT(L,\cdot)\Rightarrow\alpha=\lambda\Rightarrow\alpha\in \Lambda(L,\cdot)$.
Also, $L_{x\alpha\cdot x^\rho}=\lambda$.

$R_{x\phi^{-1}}\in R_x\mathcal{P}(L,\cdot)$ implies that $R_{x\phi^{-1}}=R_{x\rho}\Rightarrow
y\cdot x\phi^{-1}=(yx)\rho\Rightarrow (I,\phi^{-1},\rho) \in AUT(L,\cdot)\Rightarrow\phi^{-1}=\rho\Rightarrow\phi\in\mathcal{P}(L,\cdot)$.

Next, $R_{x\alpha\cdot x^\rho} \in \Phi(L,\cdot)\Rightarrow yR_{x\alpha\cdot x^\rho} =y\pi~\textrm{for some}~\pi\in\Phi(L,\cdot)\Rightarrow
y(x\alpha\cdot x^{\rho})\cdot xy=y\pi\cdot xy\Rightarrow y(x\alpha\cdot y)=y\pi\cdot xy$.
Putting $y=e$, $e\cdot (x\alpha\cdot e)=e\pi\cdot xe\Rightarrow \alpha=L_{e\pi}$.

Next, $L_{x\alpha\cdot x^\rho} \in \Psi(L,\cdot)\Rightarrow yL_{x\alpha\cdot x^\rho} =y\varrho~\textrm{for some}~\varrho\in\Psi(L,\cdot)\Rightarrow (x\alpha\cdot x^{\rho})\cdot xy=(xy)\varrho\Rightarrow  x\alpha\cdot y = (xy)\varrho$.

Thus, $(\alpha,I,\varrho) \in AUT(L,\cdot)\Rightarrow\alpha=\varrho=\lambda$.
Also, $R_{x^{\lambda}\cdot x\phi^{-1}} \in \Phi(L,\cdot)\Rightarrow R_{x^{\lambda}\cdot x\phi^{-1}} =\pi\Rightarrow
y(x^\lambda\cdot x\phi^{-1})=y\pi\Rightarrow (yx)(x^{\lambda}\cdot x\phi^{-1})=(yx)\pi( y\cdot x\phi^{-1} = (yx)\pi\Rightarrow\Rightarrow
(I,\phi^{-1},\pi) \in AUT(L,\cdot)\Rightarrow\phi^{-1}=\pi=\rho$. So, $\phi\in \Phi(L,\cdot)$.

Finally, $L_{x^{\lambda}\cdot x\phi^{-1}} \in \Psi(L,\cdot)\Rightarrow yL_{x^{\lambda}\cdot x\phi^{-1}} =y\varrho\Rightarrow
(yx)(x^\lambda\cdot x\phi^{-1})y=yx\cdot y\varrho\Rightarrow (y\cdot  x\phi^{-1})y=yx\cdot y\varrho$.

With $y=e$, then $(e\cdot x\phi^{-1})e=ex\cdot e\varrho\Rightarrow
\phi^{-1}=R_{e\varrho}\Rightarrow \phi=R_{e\varrho}^{-1}$.

Since $\alpha = \varrho =\lambda$ and $\phi^{-1}=\pi=\rho$ and $\alpha$ and
$\phi$ are arbitrary elements from $A(L)$, then $\alpha\in\Lambda(L,\cdot)$, $\alpha\in\mathcal{P}(L,\cdot)$, $\alpha\in\Phi(L,\cdot)$ and $\alpha\in\Psi(L,\cdot)$. Hence,
\begin{displaymath}
A(L)= \mathcal{P}(L,\cdot)\cap\Lambda(L,\cdot)\cap\Phi(L,\cdot)\cap\Psi(L,\cdot)
\end{displaymath}
For any $\alpha\in A(L,\cdot)$, $\alpha = L_{e\Phi}=R^{-1}_{e\Psi}$ for some
$\pi\in \Phi(L,\cdot)$ and some $\varrho\in \Psi(L,\cdot)$.

\begin{myth}\label{isere3.5}
Let $L$ be a loop with an A-holomorph Osborn loop $H(L)$. Then for all $x\in L$ and $\alpha,\phi \in A(L)$.
\begin{itemize}
\item[(a)]
$L^{-1}_xL_{x\alpha}=L_{x\alpha\cdot x^\rho}\stackrel{\delta,\beta}{\cong}x\alpha\cdot x^\rho$.
\item[(b)]
$R^{-1}_xR_{x\phi^{-1}}=L^{-1}_xR^{-1}_xR_{x\phi^{-1}}L_x=R_{x^{\lambda}\cdot x\phi^{-1}}
\stackrel{\psi,\sigma}{\cong}x^{\lambda}\cdot x\phi^{-1}=x(x^{\rho}/x \cdot x\phi^{-1})$.
\item[(c)]
$R^{-1}_xR_{x\phi^{-1}}=R_{x\alpha\cdot x^\rho}\stackrel{\sigma}{\cong}x\alpha\cdot x^\rho,L_{x^{\lambda}\cdot x\phi^{-1}}\stackrel{\beta}{\cong}x^{\lambda}\cdot x\phi^{-1}$.
\item[(d)]
$R^{-1}_xR_{x\phi^{-1}}=R_{x\alpha\cdot x^\rho}\stackrel{\varphi}{\cong}L_{x\alpha\cdot x^\rho},
R_{x^{\lambda}\cdot x\phi^{-1}}\stackrel{\varphi}{\cong}L_{x^{\lambda}\cdot x\phi^{-1}}$.
\end{itemize}
\end{myth}
{\bf Proof:}\\
Let $U=L^{-1}_xL_{x\alpha}\in\Lambda(L,\cdot)$, so $\delta(U)=eU=eL^{-1}_xL_{x\alpha}=x^{\rho}L_{x\alpha}=x\alpha\cdot x^\rho\in N_{\lambda}(L,\cdot)$. Thus, $L^{-1}_xL_{x\alpha}\stackrel{\delta}{\cong}x\alpha\cdot
 x^\rho~\forall~
x\in L, \alpha\in A(L)$.\\
Let $U= L_{x\alpha\cdot x^\rho}\in \Lambda(L,\cdot)$, then, $\delta(U)=eU=eL_{x\alpha\cdot x^\rho}=x\alpha\cdot x^\rho \in N_{\lambda}(L,\cdot)$. Thus, $L_{x\alpha\cdot x^\rho}\stackrel{\delta}{\cong} x\alpha\cdot x^\rho~\forall~x\in L,\alpha\in A(L)$.\\
Let $U=R^{-1}_xR_{x\phi^{-1}}\in \mathcal{P}(L,\cdot)$, then, $\psi(U)=eU=eR^{-1}_xR_{x\phi^{-1}}=x^{\lambda}R_{x\phi^{-1}}=x^{\lambda}\cdot x\phi^{-1}\in
N_{\rho}(L,\cdot)$. Thus, $R^{-1}_xR_{x\phi^{-1}}\stackrel{\psi}{\cong}x^\lambda\cdot x\phi^{-1}~\forall~x\in L,\phi \in A(L)$.\\
Let $U=L^{-1}_xR^{-1}_xR_{x\phi^{-1}}L_x\in \mathcal{P}(L,\cdot)$, so $\psi(U)=eU=
eL^{-1}_xR^{-1}_xR_{x\phi^{-1}}L_x=x(x^{\rho}/x \cdot x\phi^{-1})\in N_{\rho}(L,\cdot)$.
Thus, $L^{-1}_xR^{-1}_xR_{x\phi^{-1}}L_x\stackrel{\psi}{\cong}x(x^\rho/x \cdot x\phi^{-1})$.\\
Let $U=R_{x^{\lambda}\cdot x\phi^{-1}}\in \mathcal{P}(L,\cdot)$. So, $\psi(U)=eU=eR_{x^{\lambda}\cdot x\phi^{-1}}=x^{\lambda}\cdot x\phi^{-1}\in N_\rho(L)$. Thus
$R_{x^{\lambda}\cdot x\phi^{-1}}\stackrel{\psi}{\cong} {x^{\lambda}\cdot x\phi^{-1}}$.\\
Let $U=R_{x\alpha\cdot x^\rho}\in \Phi(L,\cdot)$, so  $\sigma(U)=eR_{x\alpha\cdot x^\rho}=
{x\alpha\cdot x^\rho}\in N_\mu(L,\cdot)$. Thus, $R_{x\alpha\cdot x^\rho}\stackrel{\sigma}{\cong}{x\alpha\cdot x^\rho}$.\\
Let $U=R_{x^{\lambda}\cdot x\phi^{-1}}\in \Phi(L,\cdot)$, so, $\sigma(U)=eR^{-1}_{x^{\lambda}\cdot x\phi^{-1}}={x^{\lambda}\cdot x\phi}\in N_\mu(L,\cdot)$. Thus, $R_{x^{\lambda}\cdot x\phi^{-1}}\stackrel{\sigma}{\cong}{x^{\lambda}\cdot x\phi^{-1}}$.\\
Let $U=L_{{x\alpha}\cdot x^\rho}\in \Psi(L,\cdot)$, so, $\beta(U)=eL_{{x\alpha}\cdot x^\rho}=
{x\alpha}\cdot x^\rho \in N_\mu(L,\cdot)$. Thus, $L_{{x\alpha}\cdot x^\rho}\stackrel{\beta}{\cong}{x\alpha}\cdot x^\rho$.\\
Let $U=L_{x^{\lambda}\cdot x\phi^{-1}}\in \Psi(L,\cdot)$, so, $\beta(U)=eL_{x^{\lambda}\cdot x\phi^{-1}}=x^{\lambda}\cdot x\phi^{-1}\in N_\mu(L,\cdot)$. Thus, $L_{x^{\lambda}\cdot x\phi^{-1}}\stackrel{\beta}{\cong}x^{\lambda}\cdot x\phi^{-1}$.\\
Let $U=R_{{x\alpha}\cdot x^\rho}\in \Phi(L,\cdot)$. So, $\varphi(U)=U'=L_{{x\alpha}\cdot x^\rho}$.
Thus, $R_{{x\alpha}\cdot x^\rho}\stackrel{\varphi}{\cong}L_{{x\alpha}\cdot x^\rho}\in \Psi(L,\cdot)$.\\
Let $U=R_{x^{\lambda}\cdot x\phi^{-1}}\in \Phi(L,\cdot)$. So, $\varphi(U)=U'=L_{x^{\lambda}\cdot x\phi^{-1}}$. Thus, $R_{x^{\lambda}\cdot x\phi^{-1}}\stackrel{\varphi}{\cong}L_{x^{\lambda}\cdot x\phi^{-1}}\in \Psi(L,\cdot)$.
\begin{myth}\label{isere3.6}
Let $(L,\cdot)$ be loop with an A-holomorph Osborn loop $H(L)$. Then,
\begin{itemize}
\item[(a)] \begin{equation}
\begin{diagram}
&                                                                               & x\alpha \cdot x^\rho\\
&\ruMapsto^{\sigma}_{}&\uMapsto^{\delta,\beta}_{}\\
R_{x\alpha\cdot x^\rho} &\rMapsto^{\varphi}_{\textrm{isomorphism}} & L_{x\alpha\cdot x^\rho}
\end{diagram}\in
\begin{diagram}
&                                                                               & N_\mu(L,\cdot),N_\lambda(L,\cdot)\\
&\ruTo^{\sigma}_{}&\uTo^{\delta,\beta}_{\textrm{isomorphism}}\\
\Phi(L,\cdot) &\rTo^{\varphi}_{\textrm{isomorphism}} & \Lambda(L,\cdot),\Psi(L,\cdot)
\end{diagram}
\end{equation}
for all $x\in L,\alpha\in A(L)$, i.e. $\sigma =\varphi\delta$ and $\sigma =\varphi\beta$.
\item[(b)] \begin{equation}
\begin{diagram}
&                                                                               & x^\lambda \cdot x\phi^{-1}\\
&\ruMapsto^{\psi,\sigma}_{}&\uMapsto^{\beta}_{}\\
R_{x^{\lambda}\cdot x\phi^{-1}} &\rMapsto^{\varphi}_{\textrm{isomorphism}} & L_{x^{\lambda}\cdot x\phi^{-1}}
\end{diagram}\in
\begin{diagram}
&                                                                               & N_\rho(L,\cdot),N_\mu(L,\cdot)\\
&\ruTo^{\psi,\sigma}_{}&\uTo^{\beta}_{\textrm{isomorphism}}\\
\mathcal{P}(L,\cdot),\Phi(L,\cdot) &\rTo^{\varphi}_{\textrm{isomorphism}} & \Psi(L,\cdot)
\end{diagram}
\end{equation}
for all $x\in L,\phi\in A(L)$, i.e. $\sigma =\varphi\beta$ and $\psi =\varphi\beta$.
\end{itemize}
\end{myth}
{\bf Proof:}\\
The proof follows from Theorem~\ref{isere3.5}.

\begin{myth}\label{isere3.6.1}
Let $(L,\cdot)$ be a loop with an A-holomorph Osborn loop $H(L)$.
\begin{itemize}
\item[(a)]
The commutative diagram
\begin{equation}\label{cdeq:2}
\begin{diagram}
 \Lambda(L,\cdot)&\rTo^{\delta}_{}    & N(L,\cdot)\\
\dDotsto^{\delta_1}_{}&\ruTo^{\sigma}_{}&\uTo^{\beta}_{\textrm{isomorphism}}\\
\Phi(L,\cdot) &\rTo^{\varphi}_{\textrm{isomorphism}} & \Psi(L,\cdot)
\end{diagram}
\end{equation}
is true and so, $\delta=\delta_{1}\sigma=\delta_{1}\varphi\beta$, $L_{x\alpha\cdot x^\rho}\stackrel{\delta_1}{\cong}R_{x\alpha\cdot x^\rho}$ and $x\alpha\cdot x^\rho\in Z(L,\cdot)$ for all $x\in L,\alpha\in A(L)$.
\item[(b)]
The commutative diagram
\begin{equation}\label{}
\begin{diagram}
 \Lambda(L,\cdot)&\rTo^{\delta}_{}    & N(L,\cdot)\\
\uDotsto^{\delta_2}_{}&\ruTo^{\sigma}_{}&\uTo^{\beta}_{\textrm{isomorphism}}\\
\Phi(L,\cdot) &\rTo^{\varphi}_{\textrm{isomorphism}} & \Psi(L,\cdot)
\end{diagram}
\end{equation}
is true and so, $\sigma=\delta_{2}\delta$ and $\delta_{2}\delta=\varphi\beta$ and $R_{x\alpha\cdot x^\rho}\stackrel{\delta_1}{\cong}L_{x\alpha\cdot x^\rho}$.
\item[(c)] The commutative diagram
\begin{equation}\label{cdeq:3}
\begin{diagram}
 \mathcal{P}(L,\cdot)&\rTo^{\psi}_{}    & N(L,\cdot)\\
\dDotsto^{\psi_1}_{}&\ruTo^{\sigma}_{}&\uTo^{\beta}_{\textrm{isomorphism}}\\
\Phi(L,\cdot) &\rTo^{\varphi}_{\textrm{isomorphism}} & \Psi(L,\cdot)
\end{diagram}
\end{equation}
is true and so, $\psi=\psi_{1}\sigma =\psi_{1}\varphi\beta$ and $R_{x^{\lambda}\cdot x\phi^{-1}}\stackrel{\psi_1}{\cong}R_{x^{\lambda}\cdot x\phi^{-1}}$.
\item[(d)] The commutative diagram
\begin{equation}\label{cdeq:3.1}
\begin{diagram}
 \mathcal{P}(L,\cdot)&\rTo^{\psi}_{}    & N(L,\cdot)\\
\uDotsto^{\psi_2}_{}&\ruTo^{\sigma}_{}&\uTo^{\beta}_{\textrm{isomorphism}}\\
\Phi(L,\cdot) &\rTo^{\varphi}_{\textrm{isomorphism}} & \Psi(L,\cdot)
\end{diagram}
\end{equation}
is true and so, $\sigma=\psi_{2}\psi$ and $\psi_{2}\psi =\varphi\beta$ and $R_{x^{\lambda}\cdot x\phi^{-1}}\stackrel{\psi_2}{\cong}R_{x^{\lambda}\cdot x\phi^{-1}}$.
\item[(e)] The commutative diagram
\begin{equation}\label{cdeq:5}
\begin{diagram}
 \mathcal{P}(L,\cdot)&&&&\\
&\rdDotsto~{\omega_1}\rdTo(4,2)^{\psi}\rdDotsto(2,4)_{\psi_1}&&&\\
&    &\Lambda(L,\cdot)&\rTo_{\delta}&N(L,\cdot)\\
&    &\dDotsto_{\delta_1}       &\ruTo^{\sigma}_{} &\uTo_{\beta}\\
&    &\Phi(L,\cdot) &\rTo^{\varphi}&\Psi(L,\cdot)
\end{diagram}
\end{equation}
is true and so, $\psi_1=\omega_1\delta_{1}$ and $\psi =\omega_1\delta$ and $R_{x^{\lambda}\cdot x\phi^{-1}}\stackrel{\omega_1}{\cong}L_{x^{\lambda}\cdot x\phi^{-1}}$.
\item[(f)] The commutative diagram
\begin{equation}\label{cdeq:6}
\begin{diagram}
 \mathcal{P}(L,\cdot)&&&&\\
&\luDotsto~{\omega_2}\rdTo(4,2)^{\psi}\luDotsto(2,4)_{\psi_2}&&&\\
&    &\Lambda(L,\cdot)&\rTo_{\delta}&N(L,\cdot)\\
&    &\uDotsto_{\delta_2}       &\ruTo^{\sigma}_{} &\uTo_{\beta}\\
&    &\Phi(L,\cdot) &\rTo^{\varphi}&\Psi(L,\cdot)
\end{diagram}
\end{equation}
is true and so, $\psi_2=\delta_{2}\omega_2$ and $\delta =\omega_2\psi$, and $L_{x^{\lambda}\cdot x\phi^{-1}}\stackrel{\omega_2}{\cong}R_{x^{\lambda}\cdot x\phi^{-1}}$.
\end{itemize}
\end{myth}
{\bf Proof:}\\
The proof is a consequence of Theorem~\ref{isere3.5}.
\begin{mycor}
Let $L$ be a loop with an A-holomorph Osborn loop $H(L)$.
\begin{itemize}
\item[(a)] The commutative diagram
\begin{equation}\label{}
\begin{diagram}
 \mathcal{P}(L,\cdot)&\rTo^{w_{1}=\psi\delta^{-1}}_{}    & \Lambda(L,\cdot)\\
\uTo^{\psi_1=\psi\beta^{-1}\varphi^{-1}=\psi\sigma^{-1}}_{}& &\uDotsto^{\epsilon_{1}=\beta\delta^{-1}}_{\textrm{isomorphism}}\\
\Phi(L,\cdot) &\rTo^{\varphi}_{\textrm{isomorphism}} & \Psi(L,\cdot)
\end{diagram}
\end{equation}
is true and $L_{x^{\lambda}\cdot x\phi^{-1}}\stackrel{\epsilon_1}{\cong}L_{x^{\lambda}\cdot x\phi^{-1}}$.
\item[(b)] The commutative diagram
\begin{equation}\label{}
\begin{diagram}
 \mathcal{P}(L,\cdot)&\rTo^{w_{2}=\delta\psi^{-1}}_{}    & \Lambda(L,\cdot)\\
\dTo^{\psi_2=\varphi\beta\psi^{-1}=\sigma\psi^{-1}}_{}& &\dDotsto^{\epsilon_{2}=\delta\beta^{-1}}_{\textrm{isomorphism}}\\
\Phi(L,\cdot) &\rTo^{\varphi}_{\textrm{isomorphism}} & \Psi(L,\cdot)
\end{diagram}
\end{equation}
\end{itemize}
is true and $L_{x^{\lambda}\cdot x\phi^{-1}}\stackrel{\epsilon_2}{\cong}L_{x^{\lambda}\cdot x\phi^{-1}}$.
\end{mycor}
{\bf Proof:}\\
The proof follows from Theorem~\ref{isere3.5} and Theorem~\ref{isere3.6}.
$\psi=\psi_{1}\varphi\beta \Rightarrow \psi_{1}=\psi\beta^{-1}\varphi^{-1}.\qquad
\beta=\epsilon_{1}\delta \Rightarrow \epsilon_{1}=\beta\delta^{-1}.\qquad
w_{1}=\psi_{1}\varphi\epsilon=\psi\delta^{-1}.\qquad
w_{2}\psi=\delta\Rightarrow w_{2}=\delta\psi^{-1}.\qquad
\psi_2\varphi =\varphi\beta\Rightarrow \psi_2=\varphi\beta\psi^{-1}$.

\section{Acknowledgement}
\paragraph{}
The first author wishes to express his profound gratitude and appreciation to the Management of Education
Trust Found Academic Staff Training and Development-2009(ETF AST and D)for the grant given him to
carry out this Research, as well as, to the management of Ambrose Alli University, Nigeria for their joint support
of the grant.

\end{document}